\documentclass[11pt,reqno]{amsart}

\usepackage{amsthm}
\usepackage{amscd}
\usepackage{amsfonts}
\usepackage{amssymb}
\usepackage{amsgen}
\usepackage{amsmath}
\usepackage{amsopn}
\usepackage{verbatim}
\usepackage{xypic}
\usepackage{xspace}
\usepackage{multicol}
\usepackage{url}
\usepackage{upref}

\theoremstyle{plain}
\newtheorem{thm}{Theorem}[section]
\newtheorem{lem}[thm]{Lemma}
\newtheorem{prop}[thm]{Proposition}

\newtheorem{conj}[thm]{Conjecture}

\theoremstyle{definition}

\theoremstyle{remark}
\newtheorem{rem}[thm]{Remark}

 \font\cyr=wncyr10
 \newcommand{\nc}{\newcommand}

\DeclareMathOperator{\Reg}{{Reg}}

\nc{\per}[1]{\underset{#1}{\boldsymbol \pi}\,}

 \nc{\MZV}{{\mathcal{MZV}}}

 \nc{\MT}{{\rm MT}}
  \nc{\bgz}{{\bar{\gz}}}
 \nc{\wt}{{\rm wt}}
 \nc{\wht}{{\widehat}}
 \nc{\bwg}{{\bigwedge}}
 \nc{\mmu}{{\boldsymbol{\mu}}}
 \nc{\mal}{{{\scriptstyle \maltese}}}
 \nc{\fA}{{\mathfrak A}}
 \nc{\HH}{{\mathfrak H}}
 \nc{\ra}{\rightarrow}
 \nc{\ors}{{\vec s\,}}
 \nc{\os}{{\overset}}
 \nc{\G}{{\mathbb G}}
 \nc{\Z}{{\mathbb Z}}
 \nc{\R}{{\mathbb R}}
 \nc{\N}{{\mathbb N}}
 \nc{\ZN}{{\mathbb Z_{\ge 0}}}
 \nc{\Q}{{\mathbb Q}}
 \nc{\C}{{\mathbb C}}
 \nc{\Cnn}{{\mathbb C}_{\ge 0}}
 \nc{\Cp}{{\mathbb C}_{>0}}
 \nc{\MPV}{{\mathcal{MPV}}}
 
 \nc{\tB}{{\tilde B}}
 \nc{\tI}{{\tilde I}}
 \nc{\tJ}{{\tilde J}}
 \nc{\tK}{{\tilde K}}
 \nc{\Li}{{\rm Li}}
 \nc{\suf}{{\ast\,}}
 \nc{\sufq}{{\ast_q\,}}
 \nc{\gam}{{\gamma}}
 \nc{\gG}{{\Gamma}}
 \nc{\om}{{\omega}}
 \nc{\vep}{{\varepsilon}}
 \nc{\ga}{{\alpha}}
 \nc{\gl}{{\lambda}}
 \nc{\gb}{{\beta}}
 \nc{\gd}{{\delta}}
 \nc{\gs}{{\sigma}}
 \nc{\gS}{{\Sigma}}

 \nc{\gk}{{\kappa}}
 \nc{\tgz}{{\tilde{\zeta}}}
 \nc{\gO}{{\Omega}}
 \nc{\sif}{{\mathcal S}}
 \nc{\gt}{{\tau}}
\nc{\gz}{{\zeta}}
 \nc{\Lra}{\Longrightarrow}
 \nc{\lra}{\longrightarrow}
 \nc{\fS}{{\mathfrak S}}
 \nc{\DD}{{\mathfrak D}}

 \nc{\Llra}{\Longleftrightarrow}
 \nc{\ol}{\overline}
 \nc{\lms}{\longmapsto}
 \nc{\cv}{{{\mathsf c}{\mathsf v}}}
 \nc{\zq}{{\zeta_q}}
 \nc\qup{{q\uparrow 1}}
 \nc{\us}{\underset}
 \nc{\tn}{{\tilde{n}}}
 \nc{\gD}{{\Delta}}
 \nc{\bi}{{\bf i}}
 \nc{\bfone}{{\bf 1}}
 \nc{\bfa}{{\bf a}}
 \nc{\bfb}{{\bf b}}
 \nc{\bfc}{{\bf c}}
 \nc{\bfd}{{\bf d}}
 \nc{\bfe}{{\bf e}}
 \nc{\bff}{{\bf f}}
 \nc{\bfg}{{\bf g}}
 \nc{\bfh}{{\bf h}}
 \nc{\bfi}{{\bf i}}
 \nc{\bfj}{{\bf j}}
 \nc{\bfn}{{\bf n}}
 \nc{\bfl}{{\bf l}}
 \nc{\bfk}{{\bf k}}
 \nc{\bfm}{{\bf m}}
 \nc{\bfo}{{\bf o}}
 \nc{\bfp}{{\bf p}}
 \nc{\bfq}{{\bf q}}
 \nc{\bfr}{{\bf r}}
 \nc{\tbfs}{{\tilde{\bf s}}}
 \nc{\bfs}{{\bf s}}
 \nc{\hbfs}{{\hat{\bf s}}}
 \nc{\hs}{{\hat{s}}}
 \nc{\ts}{\tilde{s}}
 \nc{\bft}{{\bf t}}
 \nc{\bfu}{{\bf u}}
 \nc{\bfv}{{\bf v}}
 \nc{\bfw}{{\bf w}}
 \nc{\bfx}{{\bf x}}
 \nc{\bfy}{{\bf y}}
 \nc{\bfz}{{\bf z}}
 \nc{\bfB}{{\bf B}}
 \nc{\bfP}{{\bf P}}
 \nc{\bfQ}{{\bf Q}}
 \nc{\bfY}{{\bf Y}}
 \nc{\bfgb}{{\boldsymbol \gb}}
 \nc{\bfga}{{\boldsymbol \ga}}
 \nc{\bfrho}{{\boldsymbol \rho}}
 \nc{\bfchi}{{\boldsymbol \chi}}
 \nc{\QX}{{\Q\langle \bfX\rangle}}
 \nc{\QY}{{\Q\langle \bfY\rangle}}
 \nc{\CX}{{\C\langle \bfX\rangle}}
 \nc{\CY}{{\C\langle \bfY\rangle}}
 \nc{\QXX}{{\Q\langle\!\langle \bfX\rangle\!\rangle}}
 \nc{\QYY}{{\Q\langle\!\langle \bfY\rangle\!\rangle}}
 \nc{\CXX}{{\C\langle\!\langle \bfX\rangle\!\rangle}}
 \nc{\CYY}{{\C\langle\!\langle \bfY\rangle\!\rangle}}

 \nc{\bbA}{{\mathbb A}}
 \nc{\bbB}{{\mathbb B}}
 \nc{\bbC}{{\mathbb C}}
 \nc{\bbD}{{\mathbb D}}
 \nc{\bbE}{{\mathbb E}}
 \nc{\bbF}{{\mathbb F}}
 \nc{\bbG}{{\mathbb G}}
 \nc{\bbH}{{\mathbb H}}
 \nc{\bbI}{{\mathbb I}}
 \nc{\bbJ}{{\mathbb J}}
 \nc{\bbK}{{\mathbb K}}
 \nc{\bbL}{{\mathbb L}}
 \nc{\bbM}{{\mathbb M}}
 \nc{\bbN}{{\mathbb N}}
 \nc{\bbO}{{\mathbb O}}
 \nc{\bbP}{{\mathbb P}}
 \nc{\bbQ}{{\mathbb Q}}
 \nc{\bbR}{{\mathbb R}}
 \nc{\bbS}{{\mathbb S}}
 \nc{\bbT}{{\mathbb T}}
 \nc{\bbU}{{\mathbb U}}
 \nc{\bbV}{{\mathbb V}}
 \nc{\bbW}{{\mathbb W}}
 \nc{\bbX}{{\mathbb X}}
 \nc{\bbY}{{\mathbb Y}}
 \nc{\bbZ}{{\mathbb Z}}
 \nc{\bba}{{\mathbb a}}
 \nc{\bbb}{{\mathbb b}}
 \nc{\bbc}{{\mathbb c}}
 \nc{\bbd}{{\mathbb d}}
 \nc{\bbe}{{\mathbb e}}
 \nc{\bbf}{{\mathbb f}}
 \nc{\bbg}{{\mathbb g}}
 \nc{\bbh}{{\mathbb h}}
 \nc{\bbi}{{\mathbb i}}
 \nc{\bbk}{{\mathbb k}}
 \nc{\bbl}{{\mathbb l}}
 \nc{\bbm}{{\mathbb m}}
 \nc{\bbn}{{\mathbb n}}
 \nc{\bbo}{{\mathbb o}}
 \nc{\bbp}{{\mathbb p}}
 \nc{\bbq}{{\mathbb q}}
 \nc{\bbr}{{\mathbb r}}
 \nc{\bbs}{{\mathbb s}}
 \nc{\bbt}{{\mathbb t}}
 \nc{\bbu}{{\mathbb u}}
 \nc{\bbv}{{\mathbb v}}
 \nc{\bbw}{{\mathbb w}}
 \nc{\bbx}{{\mathbb x}}
 \nc{\bby}{{\mathbb y}}
 \nc{\bbz}{{\mathbb z}}

 \nc{\calA}{{\mathcal A}}
 \nc{\calB}{{\mathcal B}}
 \nc{\calC}{{\mathcal C}}
 \nc{\calD}{{\mathcal D}}
 \nc{\calE}{{\mathcal E}}
 \nc{\calF}{{\mathcal F}}
 \nc{\calG}{{\mathcal G}}
 \nc{\calH}{{\mathcal H}}
 \nc{\calI}{{\mathcal I}}
 \nc{\calJ}{{\mathcal J}}
 \nc{\tcalI}{{\tilde{\mathcal I}}}
 \nc{\tcalJ}{{\tilde{\mathcal J}}}
 \nc{\calK}{{\mathcal K}}
 \nc{\calL}{{\mathcal L}}
 \nc{\calM}{{\mathcal M}}
 \nc{\calN}{{\mathcal N}}
 \nc{\calO}{{\mathcal O}}
 \nc{\calP}{{\mathcal P}}
 \nc{\calQ}{{\mathcal Q}}
 \nc{\calR}{{\mathcal R}}
 \nc{\calS}{{\mathcal S}}
 \nc{\calT}{{\mathcal T}}
 \nc{\calU}{{\mathcal U}}
 \nc{\calV}{{\mathcal V}}
 \nc{\calW}{{\mathcal W}}
 \nc{\calX}{{\mathcal X}}
 \nc{\calY}{{\mathcal Y}}
 \nc{\calZ}{{\mathcal Z}}
  \nc{\cala}{{\mathcal a}}
 \nc{\calb}{{\mathcal b}}
 \nc{\calc}{{\mathcal c}}
 \nc{\cald}{{\mathcal d}}
 \nc{\cale}{{\mathcal e}}
 \nc{\calf}{{\mathcal f}}
 \nc{\calg}{{\mathcal g}}
 \nc{\calh}{{\mathcal h}}
 \nc{\cali}{{\mathcal i}}
 \nc{\calj}{{\mathcal j}}
 \nc{\calk}{{\mathcal k}}
 \nc{\call}{{\mathcal l}}
 \nc{\calm}{{\mathcal m}}
 \nc{\caln}{{\mathcal n}}
 \nc{\calo}{{\mathcal o}}
 \nc{\calp}{{\mathsf p}}
 \nc{\calq}{{\mathcal q}}
 \nc{\calr}{{\mathcal r}}
 \nc{\cals}{{\mathcal s}}
 \nc{\calt}{{\mathcal t}}
 \nc{\calu}{{\mathcal u}}
 \nc{\calv}{{\mathcal v}}
 \nc{\calw}{{\mathcal w}}
 \nc{\calx}{{\mathcal x}}
 \nc{\caly}{{\mathcal y}}
 \nc{\calz}{{\mathcal z}}

 \nc{\frakA}{{\mathfrak A}}
 \nc{\frakB}{{\mathfrak B}}
 \nc{\frakC}{{\mathfrak C}}
 \nc{\frakD}{{\mathfrak D}}
 \nc{\frakE}{{\mathfrak E}}
 \nc{\frakF}{{\mathfrak F}}
 \nc{\frakG}{{\mathfrak G}}
 \nc{\frakH}{{\mathfrak H}}
 \nc{\frakI}{{\mathfrak I}}
 \nc{\frakJ}{{\mathfrak J}}
 \nc{\frakK}{{\mathfrak K}}
 \nc{\frakL}{{\mathfrak L}}
 \nc{\frakM}{{\mathfrak M}}
 \nc{\frakN}{{\mathfrak N}}
 \nc{\frakO}{{\mathfrak O}}
 \nc{\frakP}{{\mathfrak P}}
 \nc{\frakQ}{{\mathfrak Q}}
 \nc{\frakR}{{\mathfrak R}}
 \nc{\frakS}{{\mathfrak S}}
 \nc{\frakT}{{\mathfrak T}}
 \nc{\frakU}{{\mathfrak U}}
 \nc{\frakV}{{\mathfrak V}}
 \nc{\frakW}{{\mathfrak W}}
 \nc{\frakX}{{\mathfrak X}}
 \nc{\frakY}{{\mathfrak Y}}
 \nc{\frakZ}{{\mathfrak Z}}
 \nc{\fraka}{{\mathfrak a}}
 \nc{\frakb}{{\mathfrak b}}
 \nc{\frakc}{{\mathfrak c}}
 \nc{\frakd}{{\mathfrak d}}
 \nc{\frake}{{\mathfrak e}}
 \nc{\frakf}{{\mathfrak f}}
 \nc{\frakg}{{\mathfrak g}}
 \nc{\frakh}{{\mathfrak h}}
 \nc{\fraki}{{\mathfrak i}}
 \nc{\frakj}{{\mathfrak j}}
 \nc{\frakk}{{\mathfrak k}}
 \nc{\frakl}{{\mathfrak l}}
 \nc{\frakm}{{\mathfrak m}}
 \nc{\frakn}{{\mathfrak n}}
 \nc{\frako}{{\mathfrak o}}
 \nc{\frakp}{{\mathfrak p}}
 \nc{\frakq}{{\mathfrak q}}
 \nc{\frakr}{{\mathfrak r}}
 \nc{\fraks}{{\mathfrak s}}
 \nc{\frakt}{{\mathfrak t}}
 \nc{\fraku}{{\mathfrak u}}
 \nc{\frakv}{{\mathfrak v}}
 \nc{\frakw}{{\mathfrak w}}
 \nc{\frakx}{{\mathfrak x}}
 \nc{\fraky}{{\mathfrak y}}
 \nc{\frakz}{{\mathfrak z}}

 \nc{\sha}{{\mbox{\cyr x}}}
\nc{\slfour}{{{\mathfrak{sl}}(4)}}
 \nc{\sld}{{{\mathfrak{sl}}(d+1)}}
 \nc{\slr}{{{\mathfrak{sl}}(r+1)}}
 \nc{\slrr}{{{\mathfrak{sl}}(r+2)}}

 \nc{\uds}{{\underline{s}}}
\nc{\va}{{\vec a}}
\nc{\vb}{{\vec b}}
\nc{\vc}{{\vec c}}
\nc{\vdta}{{\vec \delta}}
\nc{\ve}{{\vec e}}
\nc{\vm}{{\vec m}}
\nc{\vp}{{\vec p}}
\nc{\vn}{{\vec n}}
\nc{\vmu}{{\vec \mu}}
\nc{\vr}{{\vec r}}
\nc{\vs}{{\vec s}}
\nc{\vt}{{\vec t}}
\nc{\vu}{{\vec u}}
\nc{\vx}{{\vec x}}
\nc{\vC}{{\vec C}}
\nc{\vv}{{\bf v}}

\begin{document}

\title[Witten multiple zeta values attached to $\slfour$]
{Witten multiple zeta values attached to $\slfour$}



\maketitle
\begin{center}
\sc{Jianqiang Zhao}\\
Department of Mathematics, Eckerd College, St. Petersburg, FL 33711, USA\\
Max-Planck Institut f\"ur Mathematik, Vivatsgasse 7, 53111 Bonn, Germany
\end{center}

\vskip0.6cm
\begin{center}
\sc{Xia Zhou}\\
Department of Mathematics, Zhejiang University, Hangzhou, P.R.\ China, 310027
\end{center}

\vskip0.6cm
\noindent{\small {\bf Abstract.}
In this paper we shall prove that every Witten multiple zeta value
of weight $w>3$ attached to $\slfour$ at nonnegative integer arguments
is a finite $\Q$-linear combination of MZVs of weight $w$ and
depth three or less, except for the nine irregular cases where
the Riemann zeta value $\gz(w-2)$ and the 
double zeta values of weight $w-1$ and depth $<3$ are also needed.
%\end{abstract}

\vskip0.6cm
\section{Introduction}
It is well-known that suitably defined zeta and $L$-functions and their
special values often play significant roles in many areas of mathematics.
In \cite{W} Witten studied one variable zeta functions
attached to various Lie algebras and related their special values
to the volumes of certain moduli spaces of vector bundles of curves.
Zagier \cite{Zag} (and independently Garoufalidis) gave direct proofs
that such functions at positive even integers are rational multiples
of powers of $\pi$. More recently Matsumoto and his
collaborators \cite{KMT0,KMT1,KMT2,KMT3,MTs}
defined multiple variable versions of these functions and
began to investigate their analytical and arithmetical properties.

Let $\N$ be the set of positive integers and $\N_0=\N\cup\{0\}$. For any
$d\in \N$ we let $[d]=(1,\dots,d)$ be a poset with the usual
increasing order. By
$\bfi=(i_1,\dots,i_\ell)\subseteq [d]$ we mean a nonempty subset of $[d]$
as a poset. We say the length of $\bfi$ is $\lg(\bfi)=\ell$
and the weight of $\bfi$ is $\wt(\bfi)=i_1+\dots+i_\ell$. We define
the \emph{generalized multiple zeta function} of depth $d$ as
\begin{equation}\label{equ:gzetaDef}
 \gz_d\big((s_\bfi)_{\bfi\subseteq [d]}\big):=
    \sum_{m_1,\dots,m_d=1}^\infty \prod_{\bfi\subseteq [d]}
    (\sum_{j=1}^{\lg(\bfi)} m_{i_j})^{-s_\bfi}.
\end{equation}
For example, the Euler-Zagier
multiple zeta function \cite{AET,Zag,Zana}
\begin{equation}\label{equ:gzDef}
  \gz(s_1,\dots,s_d):=\sum_{m_1> \dots>m_d\ge 1}
 m_1^{-s_1} m_2^{-s_2}\cdots  m_d^{-s_d}
\end{equation}
corresponds to the special case that $s_\bfi=0$ unless
$\bfi=(1,2,\dots,\ell)$ for $\ell=1,\dots,d$.
The Mordell-Tornheim multiple zeta function defined by \eqref{equ:MT}
(see \cite{Mord,Torn})
is the case where $s_\bfi=0$ unless $\lg(\bfi)=1$ or $\lg(\bfi)=d$.
If we set $s_\bfi=0$ unless $\bfi=(a,a+1,\dots,b)$ is a consecutive
string of positive integers then we get exactly Witten multiple
zeta function associated
to the special linear Lie algebra $\gz_\sld$ (see \cite{MTs}):
\begin{equation}\label{equ:sldDef}
  \gz_\sld\big((s_{i,j})_{1\le i\le j\le d}\big):=\sum_{m_1,\dots,m_d=1}^\infty
 \prod_{1\le i\le j\le d} (m_i+m_{i+1}+\cdots+m_j)^{-s_{i,j}}.
\end{equation}

The generalized multiple zeta-functions defined by \eqref{equ:gzetaDef} 
are special cases of the functions studied by Essouabri \cite{Es}, de Crisenoy \cite{Cr},
and Matsumoto \cite{Ma}. In particular we know that $\gz_d(\bfs)$
has meromorphic continuation to the whole complex space $\C^{2^d-1}$.
However, in the form \eqref{equ:gzetaDef} we may have better
control of its arithmetical properties, namely, we may be able
to compute explicitly their special values at nonnegative integers.
As usual we say a value of the
Euler-Zagier multiple zeta function at \textbf{positive} integers
a \emph{multiple zeta value} (MZV for short) if it is finite. Our
major interest is to solve the following problem.

\medskip
\noindent\textbf{Main Problem.}
Suppose $\bfs\in \N^{2^d-1}$ (resp.\ $\bfs\in \N_0^{2^d-1}$) and
$\gz_d(\bfs)$ converges. Is $\gz_d(\bfs)$ always a $\Q$-linear
combination of MZVs of the same weight (resp.\ same or lower weights)
and of depth $d$ or less?

\medskip

When $d=2$ the function in \eqref{equ:gzetaDef} becomes
Mordell-Tornheim double zeta function. By the main result of \cite{ZB} 
(see Prop.~\ref{prop:ZBred}) we know the above Main Problem has an 
affirmative answer for all Mordell-Tornheim multiple zeta functions.
In this paper we will consider $\gz_\slfour(\bfs)$ which is essentially
the case when $d=3$.

Throughout the rest of the paper whenever we say some
special value with positive (resp.\ nonnegative) integer arguments
is expressible by MZVs we mean that the value
can be expressed as a finite $\Q$-linear combination of MZVs
of the same (resp. same or lower) weights
and the same or lower depths. Our main result is
the following theorem which provides an affirmative answer
to the above Main Problem for the case $d=3$.

\begin{thm}\label{thm:depth3}
Suppose $\bfs\in \N^7$ (resp.\ $\bfs\in \N_0^7$).
If $\gz_3(\bfs)$ converges then  it is expressible by MZVs.
Moreover, very Witten multiple zeta value of weight $w$ attached
to $\slfour$ at nonnegative integers is a $\Q$-linear
combination of MZVs of weight $w$ and depth three or less,
except for the nine irregular cases defined by \eqref{equ:irr1} 
to \eqref{equ:irr5} where the Riemann zeta value $\gz(w-2)$
and the double zeta values of weight 
$w-1$ and depth $\le 2$ are also needed.
\end{thm}

The first author wishes to thank
Max-Planck-Institut f\"ur Mathematik for providing financial support
during his sabbatical leave when this work was done.
The second author is supported by the National Natural Science
Foundation of China, Project 10871169.

\section{Some preliminary results}
In this section we collect some useful facts which will be
convenient for us to present our main result in later sections.

\subsection{Convergence domain}
We assume that all components of $\bfs$ are integers and derive
the necessary and sufficient conditions for \eqref{equ:gzetaDef}
to converge although the conditions are still sufficient
if we allow complex variables and take the corresponding
real parts of the variables in the conditions.

Recall that the MZV in \eqref{equ:gzDef} converges if and only if
\begin{equation}\label{equ:criterion}
 \sum_{j=1}^\ell s_j>\ell
\end{equation}
for all $\ell=1,\dots,d$. It is straight-forward to see the same holds
for the ``star'' version of the multiple zeta function:
\begin{equation}\label{equ:gz*Def}
  \gz^*(s_1,\dots,s_d):=\sum_{m_1\ge \dots\ge m_d\ge 1}
 m_1^{-s_1} m_2^{-s_2}\cdots  m_d^{-s_d}.
\end{equation}
Special values of \eqref{equ:gz*Def} were studied in
\cite{Hoff} and \cite{OZ}.
We can extend this convergence criterion easily to our generalized
multiple zeta functions.
\begin{prop}\label{prop:conv}
The generalized MZV
\begin{equation}\label{equ:gzetaDef1}
 \gz_d(\bfs):=\sum_{m_1,\dots,m_d=1}^\infty \prod_{\bfi\subseteq [d]}
(\sum_{j=1}^{\lg(\bfi)} m_{i_j})^{-s_\bfi}
\end{equation}
converges if  and only if for all  $\ell=1,\dots,d$ and all
$\bfi=(i_1,\dots,i_\ell)\subseteq [d]$
\begin{equation}\label{equ:domainOfConv}
  \sum_{\bfj \text{ contains at least one of }i_1,\dots,i_\ell}  s_\bfj >\ell.
\end{equation}
\end{prop}
\begin{proof} The idea of the proof is similar to that of \cite[Thm.~4]{ZB}.
First we observe that for any subset $\bfi=(i_1,\dots,i_\ell)\subseteq [d]$
we have
$$\max\{m_{i_j}:1\le j\le \ell\}\le \sum_{j=1}^\ell  m_{i_j}
\le  \ell\max\{m_{i_j}:1\le j\le \ell\}.$$
Hence
\begin{equation}\label{equ:compRHS}
\sum_{\pi\in \frakS_d}\sum_{m_{\pi(1)}\ge \cdots\ge m_{\pi(d)}\ge 1}
\prod_{\bfi\subseteq [d]}
(m_{\pi(i_1)})^{-s_{\pi(\bfi)} }\le
 \sum_{m_1,\dots,m_d=1}^\infty \prod_{\bfi\subseteq [d]}
(\sum_{j=1}^{\lg(\bfi)} m_{i_j})^{-s_\bfi}
\le d \cdot \text{LHS},
\end{equation}
where LHS means the quantity at the extreme left of the above inequalities.
Observe that for each fixed $\pi\in \frakS_d$ the power of $m_{\pi(j)}$ in LHS
is $-\sum_\bfi s_{\pi(\bfi)}$ where $\bfi$ runs through all subsets of $[d]$
whose first component is $j$. Hence
the criterion \eqref{equ:criterion} implies that LHS and $d$(LHS) of
\eqref{equ:compRHS} converges if and only if for each fixed  $\pi\in \frakS_d$
\begin{equation*}
 \sum_{i_1=1}^{\ell} \  \sum_{\bfi: \text{ first component is } i_1} s_{\pi(\bfi)}>\ell ,
\end{equation*}
for all $\ell=1,\dots,d$. Let $i_j=\pi^{-1}(j)$ ($j=1,\dots,\ell$).
Then in the above sum $\bfi$ runs through all subset of $[d]$
containing at least some $i_j$ ($j=1,\dots,\ell$). This
is exactly \eqref{equ:domainOfConv}, as desired.
\end{proof}
\begin{rem} 
It is easy to see that similar result holds if we replace
the factors $\sum m_{i_j}$ in \eqref{equ:gzetaDef1}
by linear forms of $m_1,\dots,m_d$ with
nonnegative integer coefficients.
\end{rem}

\subsection{MZVs with arbitrary integer arguments}
In this paper we will mostly consider multiple zeta values with nonnegative
integer arguments as long as they converge. However,
we will prove a more general result as follows since the
inductive proof forces us to do so.
\begin{prop} Suppose $s_1,\dots,s_d\in \Z$. If
$\gz(s_1,\dots,s_d)$ converges then it
can be expressed as a $\Q$-linear combination of
MZVs (at \textbf{positive} integer arguments) of the same or
lower weights and the same or lower depths.
\end{prop}
\begin{proof} We prove the proposition by induction on
the depth. When $d=1$ we have nothing to prove.
Suppose the proposition holds for all MZVs of depth $d-1$.
Suppose further
\begin{equation}\label{equ:gzConvDom}
 s_1+\dots+s_j>j \quad\text{ for all }j=1,\dots,d,
\end{equation}
so that $\gz(s_1,\dots,s_d)$ converges.
Assume $-t=s_j\le 0$. Then by definition
$$\gz(s_1,\dots,s_d)=\sum_{m_1>\dots>m_{j-1}>m_{j+1}>\dots>m_d}
\frac{1}{m_1^{s_1} \cdots  m_{j-1}^{s_{j-1}}
m_{j+1}^{s_{j+1}} \cdots  m_d^{s_d}}
\left(\sum_{m_j=1}^{m_{j-1}-1} -\sum_{m_j=1}^{m_{j+1}} \right) m_j^t.$$
Now by the well-known formula (see, for e.g., \cite[p.~230, Thm.~1]{IR})
$$\sum_{m=1}^n m^t=\frac{1}{t+1}\Big(B_{t+1} (n+1) -B_{t+1}(0)\Big)$$
where $B_{t+1}(x)$ is the Bernoulli polynomial we immediately
see that $\gz(s_1,\dots,s_d)$ is a $\Q$-linear combination of
MZVs of the forms $\gz(s_1,\dots,s_{j-1}-s,s_{j+1},\dots,s_d)$ and
$\gz(s_1,\dots,s_{j-1}-s,s_{j+1},\dots,s_d)$ where
$s=0,1,\dots,-s_j+1$. All of these MZVs are easily to be shown as
convergent values by \eqref{equ:gzConvDom} and therefore the
proposition follows from induction assumption.
\end{proof}

\subsection{Mordell-Tornheim zeta functions}\label{sec:MT}
They are defined by (see \cite{Mord,Torn})
\begin{equation}\label{equ:MT}
\gz_\MT(s_1,\dots,s_d;s):=\sum_{m_1,\dots,m_d=1}^\infty
 m_1^{-s_1} m_2^{-s_2}\cdots  m_d^{-s_d}(m_1+\cdots +m_d)^{-s}.
\end{equation}
The main result of \cite{ZB} is the following
\begin{prop} \label{prop:ZBred} \emph{(\cite[Thm.~5]{ZB})}
Let $s_1,\dots,s_d$ and $s$ be nonnegative integers. If at most
one of them is equal to $0$ then the Mordell-Tornheim zeta value
$\gz_\MT(s_1,\dots,s_d;s)$ can be expressed as a $\Q$-linear
combination of MZVs of the same weight and depth.
\end{prop}
In this paper we will only need this proposition when the depth is three.

\subsection{A combinatorial lemma}
The next lemma will be used heavily throughout the paper.
\begin{lem} \label{lem:combLem}\emph{ (\cite[Lemma 1]{ZB})}
Let $r$ and $n_1,\dots,n_r$ be positive integers, and let
$x_1,\dots,x_r$ be non-zero real numbers such that $x_1+\dots+x_r\ne 0$.
Then
$$\prod_{j=1}^r \frac{1}{x_j^{n_j}} =
\sum_{j=1}^r\Bigg(\prod_{\substack{k=1\\ k\ne j}}^r
\sum_{a_k=0}^{n_k-1}\Bigg) \frac{M_j}{x^{n_j+A_j}}
    \prod_{\substack{k=1\\ k\ne j}}^r \frac{1}{x_k^{n_k-a_k}},$$
where the multi-nomial coefficient
$$M_j=\frac{(n_j+A_j-1)!}{(n_j-1)!} \prod_{\substack{k=1\\ k\ne j}}^r \frac{1}{a_k!} \qquad \text{and}\qquad A_j=\sum_{\substack{k=1\\ k\ne j}}^r a_k.$$
The notation $\displaystyle \prod_{\substack{k=1\\ k\ne j}}^r
\sum_{a_k=0}^{n_k-1}$ means the multiple sum
$\displaystyle \sum_{a_1=0}^{n_1-1}\dots\sum_{a_{j-1}=0}^{n_{j-1}-1}
\sum_{a_{j+1}=0}^{n_{j+1}-1}\dots\sum_{a_r=0}^{n_r-1}.$
\end{lem}

\section{Proof of Theorem \ref{thm:depth3}}
By definition \eqref{equ:gzetaDef}
\begin{equation}\label{equ:gzeta3Def}
 \gz_3(s_1,\dots,s_7):=\sum_{m_1,\dots,m_3=1}^\infty
 \frac{1}{m_1^{s_1}m_2^{s_2}m_3^{s_3}(m_1+m_2)^{s_4}
 (m_2+m_3)^{s_5}(m_1+m_3)^{s_6}(m_1+m_2+m_3)^{s_7}}.
\end{equation}
To guarantee convergence by \eqref{equ:domainOfConv} we need to assume:
\begin{equation}\label{equ:convConditionOfDepth3}
\left\{
\aligned
 \bfi=(1):&\quad s_1+s_4+s_6+s_7 > 1,\\
 \bfi=(1):&\quad s_2+s_4+s_5+s_7 > 1,\\
 \bfi=(1):&\quad s_3+s_5+s_6+s_7 > 1,\\
 \bfi=(1,2):&\quad s_1+s_2+s_4+s_5+s_6+s_7 > 2,\\
 \bfi=(1,3):&\quad s_1+s_3+s_4+s_5+s_6+s_7 > 2,\\
 \bfi=(2,3):&\quad s_2+s_3+s_4+s_5+s_6+s_7 > 2,\\
 \bfi=(1,2,3):&\quad s_1+s_2+s_3+s_4+s_5+s_6+s_7 >3.
\endaligned\right.
\end{equation}

We now use a series of reductions to prove the theorem. All of the
steps will be explicitly given so that one may carry out the computation
of \eqref{equ:gzeta3Def} by following them.

\medskip
\noindent
\underline{Step (i)}.
Since $2(m_1+m_2+m_3)=(m_1+m_2)+(m_2+m_3)+(m_1+m_3)$
we see clearly that we can assume either $s_4=0$ or
$s_5=0$ or $s_6=0$. In fact in Lemma~\ref{lem:combLem}
taking $r=3$, $x_1=m_1+m_2$, $x_2=m_2+m_3$, $x_3=m_1+m_3$,
$n_1=s_4$, $n_2=s_5$ and $n_3=s_6$ we get:
 \begin{align*}
\  \gz_3(s_1,\dots,s_7)=&\sum_{a_4=0}^{s_4-1}  \sum_{a_5=0}^{s_5-1}
\frac{(s_6+a_4+a_5-1)!}{2^{s_7+a_4+a_5}(s_6-1)!a_4!a_5!}
\gz_3(s_1,s_2,s_3,s_4-a_4,s_5-a_5,0,s_6+s_7+a_4+a_5)\\
+&\sum_{a_4=0}^{s_4-1}  \sum_{a_6=0}^{s_6-1}
\frac{(s_5+a_4+a_6-1)!}{2^{s_7+a_4+a_6}(s_5-1)!a_4!a_6!}
\gz_3(s_1,s_2,s_3,s_4-a_4,0,s_6-a_6,s_5+s_7+a_4+a_6)\\
+&\sum_{a_5=0}^{s_5-1}  \sum_{a_6=0}^{s_6-1}
\frac{(s_4+a_5+a_6-1)!}{2^{s_7+a_5+a_6}(s_4-1)!a_5!a_6!}
\gz_3(s_1,s_2,s_3,0,s_5-a_5,s_6-a_6,s_4+s_7+a_5+a_6).
 \end{align*}
We recommend the interested reader to check the convergence
of the above values by \eqref{equ:convConditionOfDepth3}.
The rule of thumb is that if we apply Lemma~\ref{lem:combLem}
with each $x_j$ a positive combination of indices then
the convergence is automatically guaranteed.
In each of the following steps we often omit this
convergence checking since it is straight-forward
in most cases. The only exception is \eqref{equ:limit1}
which in fact poses the most difficulty.

By symmetry we see that without loss of generality
we only need to show that
$\gz_3(s_1,\dots,s_5,0,s_6)$
is expressible by MZVs. This is nothing but the
Matsumoto's version of Witten multiple zeta function of depth 3
associated to the
special linear Lie algebra $\slfour$ (see \eqref{equ:sldDef}):
$$\gz_\slfour(s_1,\dots,s_6):=
\sum_{m_1,\dots,m_3=1}^\infty
 \frac{1}{m_1^{s_1}m_2^{s_2}m_3^{s_3}(m_1+m_2)^{s_4}
 (m_2+m_3)^{s_5}(m_1+m_2+m_3)^{s_6}}.$$
 
Before going on we need to define the so called \emph{regular} and \emph{irregular}
special values of $\gz_\slfour(s_1,\dots,s_6)$. 
Clearly the following special values are expressible by MZVs of 
mixed weights: ($\{0\}_k$ means to repeat $0$ $k$ times)
\begin{align}\label{equ:irr1}
& \gz_\slfour(\{0\}_3,b,0,a)=\gz_\slfour(\{0\}_4,b,a)= \gz(a,b,0)=
\left\{ \begin{array}{ll}
  \gz(a,b-1)-\gz(a,b), & a>1,b> 1,\\
  \gz(a-1)-\gz(a)-\gz(a,1), &  a>1,b=1,\\
\end{array} \right. \\
& \gz_\slfour(b,\{0\}_4,a)= \gz_\slfour(0,b,\{0\}_3,a)= 
\gz_\slfour(0,0,b,0,0,a)= \gz(a,0,b) \notag\\
& \ \hskip5.17cm =
\left\{ \begin{array}{ll}
 \gz(a-1,b)-\gz(a,b-1)-\gz(a,b),  & a>2,b>1,\\
 \gz(a-1,1)-\gz(a-1)+\gz(a)-\gz(a,1),  & a>2,b=1,\\
 \end{array} \right. \label{equ:irr2} \\
&  \gz_\slfour(\{0\}_5,a)=\gz(a,0,0)=\frac12\gz(a-2)-\frac32\gz(a-1)+\gz(a),
    \hskip3.5cm  a>3.\label{equ:irr3}
\end{align}}
By Step (ii) we will see that special values
\begin{equation}  \label{equ:irr4}
\gz_\slfour(0,0,s_3,s_4,0,s_6)\qquad\text{and} \qquad\gz_\slfour(s_1,0,0,0,s_5,s_6) 
\end{equation}
are also expressible by MZVs of mixed weights. Further, if $a,b\ge 2$ then we have
\begin{align} \label{equ:irr5}
 \gz_\slfour(0,0,0,a,b,0)=&\sum_{m_1,m_2,m_3=1}^\infty
  \frac{1}{(m_1+m_2)^a (m_2+m_3)^b}\\
 =& \sum_{m_2=1}^\infty\left( \sum_{m_1<m_3}
 + \sum_{m_1=m_3}+ \sum_{m_1>m_3} \right)
 \frac{1}{(m_1+m_2)^a (m_2+m_3)^b} \notag\\
 =& \gz(a,b,0)+\gz(b,a,0)+\gz(a+b,0) \notag\\
 =&\gz(b,a-1)+\gz(a,b-1)+\gz(a+b-1)-\gz(a)\gz(b).\notag
\end{align}
We call the values appearing in
the nine cases from \eqref{equ:irr1} to \eqref{equ:irr5} \emph{irregular} values. 
Otherwise $\gz_\slfour(s_1,\dots,s_6)$ is called a \emph{regular} value.
 
\medskip
\noindent
\underline{Step (ii)}. Taking $x_1=m_1$ and $x_2=m_2+m_3$
in Lemma \ref{lem:combLem} we get
\begin{equation}\label{equ:step2}
\aligned
 \gz_\slfour(s_1,\dots,s_6)
=&\sum_{a_5=0}^{s_5-1}
{s_1+a_5-1\choose a_5}
\gz_\slfour(0,s_2,s_3,s_4,s_5-a_5,s_6+s_1+a_5)\\
+&\sum_{a_1=0}^{s_1-1}
{s_5+a_1-1\choose a_1}
\gz_\slfour(s_1-a_1,s_2,s_3,s_4,0,s_6+s_5+a_1).
\endaligned
\end{equation}
We see clearly that we can assume either (ii.1): $s_5=0$ or
(ii.2): $s_1=0$ and $s_5\ge 1$.
Moreover, by taking $s_2=s_3=s_4=0$
in \eqref{equ:step2} we see that  $\gz_\slfour(s_1,0,0,0,s_5,s_6)$ in \eqref{equ:irr4} 
can be expressed by $\Q$-linear combinations of MZVs appeared in the 
\eqref{equ:irr1} and \eqref{equ:irr2}.
The argument is similar for $\gz_\slfour(0,0,s_3,s_4,0,s_6)$ in \eqref{equ:irr4}. 
Therefore in what follows we assume that $\gz_\slfour(s_1,\dots,s_6)$ 
are always regular and show that they can be expressed by $\Q$-linear
combinations of MZVs of the same weight and depth three or less.

\medskip
\noindent
\underline{Step (ii.1)}. Let $s_5=0$. Then we must have either $s_1\ge 1$ or $s_2\ge 1$
since we assume
$$\gz_\slfour(s_1,\dots,s_4,0,s_6)
=\sum_{m_1,\dots,m_3=1}^\infty
 \frac{1}{m_1^{s_1}m_2^{s_2}m_3^{s_3}(m_1+m_2)^{s_4}
 (m_1+m_2+m_3)^{s_6}}$$
is regular. Then we may use $x_1=m_1$ and $x_2=m_2$ in Lemma \ref{lem:combLem} to get
 \begin{align*}
 \gz_\slfour(s_1,\dots,s_4,0,s_6)
=&
\sum_{a_1=0}^{s_1-1}
{s_2+a_1-1\choose a_1}
\gz_\slfour(s_1-a_1,0,s_3,s_2+s_4+a_1,0,s_6)\\
+&\sum_{a_2=0}^{s_2-1}
{s_1+a_2-1\choose a_2}
\gz_\slfour(0,s_2-a_2,s_3,s_1+s_4+a_2,0,s_6) .
 \end{align*}
By symmetry we only need to consider
\begin{equation} 
\label{equ:lastStepin211} 
\gz_\slfour(0,s_2,s_3,s_4,0,s_6)
=\sum_{m_1,\dots,m_3=1}^\infty
 \frac{1}{ m_2^{s_2}m_3^{s_3}(m_1+m_2)^{s_4}
 (m_1+m_2+m_3)^{s_6}}
\end{equation}
where $s_2\ge 1$.
But now we may take $x_1=m_1+m_2$ and $x_2=m_3$
in Lemma \ref{lem:combLem} to reduce it to
either $s_3=0,s_4\ge 1$ or $s_4=0,s_3\ge 1$. Then we 
get the following two kind of values:
$$\gz_\slfour(0,s_2,0,s_4,0,s_6)=\zeta(s_6,s_4,s_2),\quad
 \gz_\slfour(0,s_2,s_3,0,0,s_6)=\gz_\MT(s_2,s_4,0;s_6),$$
where $s_2,s_3,s_4,s_6\ge 1$. But 
$\gz_\MT(s_2,s_4,0;s_6)$ is expressible by MZVs of the same weight 
and of depth three by Prop.~\ref{prop:ZBred}. Case (ii.1) is proved.

\medskip
\noindent
\underline{Step (ii.2)}. Let $s_1=0$ and $s_5\ge 1$. Consider
 $$\gz_\slfour(0,s_2,\dots,s_6)=
\sum_{m_1,\dots,m_3=1}^\infty
 \frac{1}{m_2^{s_2}m_3^{s_3}(m_1+m_2)^{s_4}
 (m_2+m_3)^{s_5} (m_1+m_2+m_3)^{s_6}}.$$
To guarantee convergence we must have $s_4+s_6>1$, $s_3+s_5+s_6>1$, 
$s_2+s_4+s_5+s_6>2$, $s_3+s_4+s_5+s_6>2$, and $s_2+s_3+s_4+s_5+s_6>3$.
Taking  $x_1=m_1+m_2$ and $x_2=m_3$ in Lemma \ref{lem:combLem} we get
\begin{align}
 \gz_\slfour(0,s_2,\dots,s_6)
=&
\sum_{a_3=0}^{s_3-1}
{s_4+a_3-1\choose a_3}
\gz_\slfour(0,s_2,s_3-a_3,0,s_5,s_4+s_6+a_3) \label{equ:case2line1}\\
+&\sum_{a_4=0}^{s_4-1}
{s_3+a_4-1\choose a_4}
\gz_\slfour(0,s_2,0,s_4-a_4,s_5,s_3+s_6+a_4) .\label{equ:case2line2}
 \end{align}
So we may assume that either (ii.2.1): $s_1=s_4=0,s_5,s_6\ge 1$ or
(ii.2.2): $s_1=s_3=0$ and $s_2,s_4,s_5\ge 1$, or
(ii.2.3): $s_1=s_2=s_3=0$ and $s_4,s_5\ge 1$. Here (ii.2.1) comes from
\eqref{equ:case2line1} while (ii.2.2) and (ii.2.3) come from
\eqref{equ:case2line2}.

\medskip
\noindent
\underline{Step (ii.2.1)}. Let $s_1=s_4=0,s_5,s_6\ge 1$. Then we must
have either $s_2\ge 1$ or $s_3\ge 1$ in
$$\gz_\slfour(0,s_2,s_3,0,s_5,s_6)=
\sum_{m_1,\dots,m_3=1}^\infty
 \frac{1}{m_2^{s_2}m_3^{s_3}
 (m_2+m_3)^{s_5} (m_1+m_2+m_3)^{s_6}}.$$
since we consider only regular values only.
Thus we may put $x_1=m_2$ and $x_2=m_3$ in Lemma \ref{lem:combLem}
to get
\begin{align*}
\gz_\slfour(0,s_2,s_3,0,s_5,s_6)
=&
\sum_{a_2=0}^{s_2-1}
{s_3+a_2-1\choose a_2}
\gz_\slfour(0,s_2-a_2,0,0,s_3+s_5+a_2,s_6)\\
+&\sum_{a_3=0}^{s_3-1}
{s_2+a_3-1\choose a_3}
\gz_\slfour(0,0,s_3-a_3,0,s_2+s_5+a_3,s_6)\\
=&
\sum_{a_2=0}^{s_2-1}
{s_3+a_2-1\choose a_2}
\gz(s_6,s_3+s_5+a_2,s_2-a_2)\\
+&\sum_{a_3=0}^{s_3-1}
{s_2+a_3-1\choose a_3}
\gz(s_6,s_2+s_5+a_3,s_3-a_3).
 \end{align*}
It is easy to see that all the triple zeta values above have the same weight.
We remind the reader that to determine the weight of a MZV it is not
enough just to add up all the components. One also needs to check that 
every component is positive.

\medskip
\noindent
\underline{Step (ii.2.2)}. Let $s_1=s_3=0$ and $s_2, s_4,s_5\ge 1$. 
By \eqref{equ:convConditionOfDepth3}, to guarantee convergence of 
\begin{equation*}
 \gz_\slfour(0,s_2,0,s_4,s_5,s_6)=
\sum_{m_1,\dots,m_3=1}^\infty
 \frac{(-1)^{s_2} }{(-m_2)^{s_2} (m_1+m_2)^{s_4}
 (m_2+m_3)^{s_5} (m_1+m_2+m_3)^{s_6}}
\end{equation*}
we need to assume 
\begin{equation}\label{equ:s1s3=0}
s_4+s_6>1,\ s_5+s_6>1,\ s_4+s_5+s_6>2,\ s_2+s_4+s_5+s_6>3.
\end{equation}
Moreover, since $\gz_\slfour(0,s_2,0,s_4,s_5,s_6)$ is regular we must
have $s_2+s_6\ge 1$.
Putting $x_1=-m_2$, $x_2=m_1+m_2$ and $x_3=m_2+m_3$
in Lemma \ref{lem:combLem} we get
 \begin{align}
\ & \gz_\slfour(0,s_2,0,s_4,s_5,s_6) \label{equ:s1s3=0Red}\\
=&\sum_{a_4=0}^{s_4-1}  \sum_{a_5=0}^{s_5-1}
\frac{(s_2+a_4+a_5-1)!}{(s_2-1)!a_4!a_5!}(-1)^{s_2}
\gz_\slfour(0,0,0,s_4-a_4,s_5-a_5,s_2+s_6+a_4+a_5)\notag\\
+&\sum_{a_2=0}^{s_2-1}  \sum_{a_5=0}^{s_5-1}
\frac{(s_4+a_2+a_5-1)!}{(s_4-1)!a_2!a_5!}(-1)^{a_2}
\gz_\slfour(0,s_2-a_2,0,0,s_5-a_5,s_4+s_6+a_2+a_5)\notag\\
+&\sum_{a_2=0}^{s_2-1}  \sum_{a_4=0}^{s_4-1}
\frac{(s_5+a_2+a_4-1)!}{(s_5-1)!a_2!a_4!}(-1)^{a_2}
\gz_\slfour(0,s_2-a_2,0,s_4-a_4,0,s_5+s_6+a_2+a_4).\notag
\end{align}
We point out that since we have used $x_1=-m_2$ we need to
check the convergence of all the above three kinds of values
under the assumption \eqref{equ:s1s3=0} even though
the checking itself is trivial.

Returning to the reduction of \eqref{equ:s1s3=0Red} we see that
the last two sums are expressible by MZVs
so we only need to consider those values appearing in the first sum,
namely, those of the form in the next case.

\medskip
\noindent
\underline{Step (ii.2.3)}. Let $s_1=s_2=s_3=0$ and $s_4,s_5\ge 1$. Then we
must have $s_6\ge 1$ since we only consider regular values. To 
guarantee convergence of
\begin{equation}\label{equ:s4s5s6}
 \gz_\slfour(0,0,0,s_4,s_5,s_6)
=\sum_{m_1,\dots,m_3=1}^\infty
 \frac{1}{ (m_1+m_2)^{s_4} (m_2+m_3)^{s_5}(m_1+m_2+m_3)^{s_6}}
\end{equation}
we need to assume $s_4+s_6>1,$ $s_5+s_6>1,$ and $s_4+s_5+s_6>3.$
If either $s_4>1$ or $s_5>1$ then we may further assume that
$s_4>1$ by change of index $m_1\leftrightarrow m_3$ in \eqref{equ:s4s5s6}.
We may set $x_1=-(m_2+m_3)$ and $x_2=m_1+m_2+m_3$ in Lemma \ref{lem:combLem}
to get
\begin{align}
\gz_\slfour(0,0,0,s_4,s_5,s_6)
&=
\sum_{a_5=0}^{s_5-2}
{s_6+a_5-1\choose a_5}(-1)^{a_5}
\gz_\slfour(s_6+a_5,0,0,s_4,s_5-a_5,0)\label{equ:s6>0} \\
&+\sum_{a_6=0}^{s_6-2}
{s_5+a_6-1\choose a_6} (-1)^{s_5}
\gz_\slfour(s_5+a_6,0,0,s_4,0,s_6-a_6) \notag\\
& \aligned
+{s_6+s_5-2\choose s_5-1}(-1)^{s_5}
\lim_{N\to \infty}
& \Big(\gz^{(N)}_\slfour(s_6+s_5-1,0,0,s_4,0,1)  \\
&  -\gz^{(N)}_\slfour(s_6+s_5-1,0,0,s_4,1,0)\Big),
\endaligned \label{equ:limit1}
\end{align}
where $\gz^{(N)}$ is the partial sum of \eqref{equ:gzetaDef} when
each index $m_{i_j}$ goes from $1$ to $N$. Observe that all the values
in the first two lines above are convergent and the second line
is already expressible by MZVs of the same weight. 
So we are reduced to consider the following two kinds of values:
\begin{enumerate}
  \item[(A).] $\gz_\slfour(s_1,0,0,s_4,t,0)$
    (with $s_1,s_4\ge 1$ and $t>1$) from first line above and

  \item[(B).] The limit \eqref{equ:limit1} (with $s:=s_5+s_6-1\ge 1$).
\end{enumerate}
We now need to divide into two subcases to compute
\eqref{equ:s6>0} (leading to (A)) and \eqref{equ:limit1} (leading to (B)):
(ii.2.3.1) $s_1=s_2=s_3=0$, either $s_4>1$ or $s_5>1$,
and (ii.2.3.2) $s_1=s_2=s_3=0$ and $s_4=s_5=1$.
Before treating the two cases separately we first deform the expression
inside the limit \eqref{equ:limit1} as follows:
\begin{align}
\ &\Big(\gz^{(N)}_\slfour(s,0,0,s_4,0,1)
 -\gz^{(N)}_\slfour(s,0,0,s_4,1,0)\Big)\notag\\
=&\sum_{m_1,m_2=1}^{N} \left(\sum_{n=1}^{m_1+m_2+N}-\sum_{n=1}^{m_1+m_2}\right)
 \frac{1}{m_1^s(m_1+m_2)^{s_4} n}-\sum_{m_1,m_2=1}^{N} \left(\sum_{n=1}^{m_2+N}-\sum_{n=1}^{m_2}\right)
 \frac{1}{m_1^s(m_1+m_2)^{s_4} n}\notag\\
=&-\sum_{m_1,m_2=1}^{N} \sum_{n=1+m_2}^{m_1+m_2}
 \frac{1}{m_1^s(m_1+m_2)^{s_4} n}
 +\sum_{m_1,m_2=1}^{N} \sum_{n=1+m_2+N}^{m_1+m_2+N}
 \frac{1}{m_1^s(m_1+m_2)^{s_4} n}. \label{equ:to0}
\end{align}

\medskip
\noindent
\underline{Step (ii.2.3.1)}. Let $s_1=s_2=s_3=0$, $s_6\ge 1$, either $s_4>1$ or $s_5>1$.
Without loss of generality assume $s_4>1$
as we have remarked after \eqref{equ:s4s5s6}. 

(A). Since $s_1\ge 1, s_4,t>1$ we have
\begin{align*}
 \gz_\slfour(s_1,0,0,s_4,t,0)=&
\sum_{m_1,m_2,m_3=1}^\infty
 \frac{1}{m_1^{s_1} (m_1+m_2)^{s_4} (m_2+m_3)^t}\\
  =&\sum_{m_1,m_2=1}^\infty\sum_{n>m_2}
  \frac{1}{m_1^{s_1} (m_1+m_2)^{s_4} n^t}\\
 =&\gz(t)\gz(s_4,s_1)-\gz_\MT(s_1,t;s_4) -\sum_{m_1,m_2=1}^\infty\sum_{n<m_2}
  \frac{1}{m_1^{s_1} (m_1+m_2)^{s_4} n^t}\\
 =&\gz(t)\gz(s_4,s_1)-\gz_\MT(s_1,t;s_4)
 -\sum_{m_1,n,m_3=1}^\infty
  \frac{1}{m_1^{s_1} (m_1+n+m_3)^{s_4} n^t}
\end{align*}
by setting $m_2=n+m_3$. But the last sum is just
$\gz_\MT(s_1,t,0;s_4)$ which is expressible by MZVs of
the same weight and of depth three by Prop.~\ref{prop:ZBred}.  

(B). The limit \eqref{equ:limit1}. This has been reduced to \eqref{equ:to0}
with $s\ge 1$ and $s_4\ge 2$. Then the second sum in \eqref{equ:to0} is bounded by
$$\frac{1}{1+N} \sum_{m_1,m_2=1}^N \frac{1}{m_1m_2}
=O(\log^2N/N)\to 0, \quad \text{as }N\to \infty.$$
By setting $m_2=m_3+n$ in
the first sum in \eqref{equ:to0} we get
\begin{align}
\ & \lim_{N\to \infty}
\Big(\gz^{(N)}_\slfour(s,0,0,{s_4},0,1)
 -\gz^{(N)}_\slfour(s,0,0,{s_4},1,0)\Big)\notag \\
=&-\sum_{m_1}^\infty  \sum_{n=1}^\infty
\left(\sum_{n=1}^{m_1+m_2}-\sum_{n=1}^{m_2} \right)
 \frac{1}{m_1^s (m_1+m_2)^{s_4} n}\notag\\
=&-\sum_{m_1,m_2=1}^\infty
\left(\sum_{n=1}^{m_1-1}+\sum_{n=m_1}+\sum_{m_1<n<m_1+m_2}+\sum_{n=m_1+m_2}
-\sum_{n=1}^{m_2} \right)
 \frac{1}{m_1^s (m_1+m_2)^{s_4} n} \label{equ:sumbreak}\\
=&\gz_\MT(s,1;{s_4})-\big(\gz({s_4},s,1)+\gz({s_4},s+1)+\gz({s_4},1,s)+\gz({s_4}+1,s)\big)
 +\sum_{m_1,m_2=1}^\infty \sum_{n=1}^{m_2-1}
 \frac{m_1^{-s} n^{-1}} {(m_1+m_2)^{s_4} }.\notag
\end{align}
In the last sum setting $m_2=n+m_3$ then we get
$$\sum_{m_1,m_2=1}^\infty \sum_{n=1}^{m_2-1}
 \frac{1}{m_1^s (m_1+m_2)^{s_4} n}=\sum_{n,m_1,m_3=1}^\infty
 \frac{1}{m_1^s (m_1+m_3+n)^{s_4} n}=\gz_\MT(s,1,0;{s_4})$$
which is expressible by MZVs of
the same weight and of depth three by Prop.~\ref{prop:ZBred}.

\medskip
\noindent
\underline{Step (ii.2.3.2)}. Let $s_1=s_2=s_3=0$ and $s_4=s_5=1$.
Then we must have $s_6\ge 2$ in order to have convergent values. 
Moreover, \eqref{equ:s6>0} is vacuous so we don't need to consider (A).
 
(B). Observe that the second term in \eqref{equ:to0} is bounded by
$$  \frac{1}{1+N} \sum_{m_1,m_2=1}^{N}\frac{1}{m_1^{s-1}m_2}=O(\log^2(N)/N)\to 0, \quad \text{as }N\to \infty$$
since $s\ge 2.$ So
$$ \lim_{N\to \infty}
\Big(\gz^{(N)}_\slfour(s,0,0,1,0,1)
 -\gz^{(N)}_\slfour(s,0,0,1,1,0)\Big)
 =-\lim_{N\to \infty} \sum_{m_1,m_2=1}^{N} \sum_{n=1+m_2}^{m_1+m_2}
 \frac{1}{m_1^s (m_1+m_2) n}$$
is expressible by MVZs by Lemma \ref{lem:techLemma}
(with $t=1$) in \S\ref{sec:reg}.

This finishes the proof of the theorem. \hfill$\Box$

\section{Regularized MZVs and a technical lemma}\label{sec:reg}
We need to derive several formulas concerning MVZs. They are
best understood as consequences of the stuffle (also called
quasi-shuffle or harmonic product) relations for the MZVs
including the divergent ones suitably regularized.

By \cite{IKZ} one may define a regularized value $\bgz(\bfs)$ (denoted
by $Z^\sha_{\bfs}(T)$ in \cite{IKZ})
if $\bfs=(s_1,\dots,s_d)$ with $s_1=1$ such that the
following stuffle relation holds:
\begin{equation}\label{equ:stuffle}
 \bgz(\bfs_1)\bgz(\bfs_2)=\bgz(\bfs_1*\bfs_2)
\end{equation}
where the the product $*$ is the stuffle product (see \cite{Hoff}
or \cite[\S1]{IKZ}).
\begin{lem} \label{lem:EulerDec}
Let $s$ and $t$ be two positive integers greater than 1. Then
\begin{equation}\label{equ:EulerDec}
 \gz(s)\gz(t)=\sum_{a=0}^s {t-1+a\choose a} \gz(t+a,s-a)+
 \sum_{b=0}^t {s-1+b\choose b} \gz(s+b,t-b).
\end{equation}
Further
\begin{equation}\label{equ:EulerDec=1}
 \gz(s+1)=\sum_{a=1}^{s-1}  \gz(1+a,s-a) .
\end{equation}
\end{lem}
\begin{proof} Equation \eqref{equ:EulerDec} is the famous
Euler's decomposition formula \cite{Eu} (or see \cite{BBBL}).
To derive \eqref{equ:EulerDec=1} we use the regularized MZVs as follows.
For any positive integer $N$ we have by Lemma~\ref{lem:combLem}
$$\sum_{m=1}^N \sum_{n=1}^{N-m} \frac{1}{mn^s}
=\sum_{m=1}^N \sum_{n=1}^{N-m}
\left(\frac{1}{m(m+n)^s}+
\sum_{a=0}^{s-1}\frac{1}{n^{s-a}(m+n)^{a+1}}\right).$$
Taking the regularization we get
$$\bgz(1)\gz(s)=\bgz(1,s)+\gz(s,1)+\sum_{a=1}^{s-1} \gz(1+a,s-a).$$
On the other hand by stuffle \eqref{equ:stuffle}
$$\bgz(1)\gz(s)=\bgz(1,s)+\gz(s,1)+\gz(s+1).$$
Equation \eqref{equ:EulerDec=1} follows immediately.
\end{proof}

\begin{lem}\label{lem:techLemma} For all positive integers $s$
and $t$ such that $s+t\ge 3$ the value
\begin{equation}\label{equ:techLemma}
 \sum_{m_1,m_2=1}^\infty \sum_{n=1+m_2}^{m_1+m_2}
 \frac{1}{m_1^s (m_1+m_2) n^t}
\end{equation}
is expressible by MVZs of weight $s+t+1$ and depth three or less.
\end{lem}
\begin{proof}
Let's consider the partial sum
$$ \sum_{1\le m_1<m_3\le N} \sum_{n=1+m_3-m_1}^{m_3}
 \frac{1}{m_1^s m_3 n^t}=\sum_{1\le m_1<m_3\le N} \left( \sum_{n=1}^{m_3}-\sum_{n=1}^{m_3-m_1}\right) \frac{1}{m_1^s m_3^t n}$$
whose limit when $N\to \infty$ is clearly \eqref{equ:techLemma}.
Using regularization and playing similar trick of breaking the sum
as we did in \eqref{equ:sumbreak} we find that
\begin{align}
\eqref{equ:techLemma}=& \Reg\Bigg\{ \sum_{1\le m_1<m_3\le N} \left( \sum_{n=1}^{m_3}-\sum_{n=1}^{m_3-m_1}\right) \frac{1}{m_1^s m_3^t n} \Bigg\}, \notag\\
 =& \bgz(1,t,s)+\bgz(1,s+t)+\bgz(1,s,t)+\gz(t+1,s)
 -\Reg\Bigg\{\sum_{1\le m_1<m_3\le N} \sum_{n=1}^{m_3-m_1}
 \frac{1}{m_1^s m_3^t n} \Bigg\}.\label{equ:techLemmaAllgzs}
\end{align}
Now the untreated sum can be deformed as follows: taking $m_2=m_3-m_1$
we get
 \begin{multline}\label{equ:reg1}
  \sum_{1\le m_1<m_1+m_2\le N} \sum_{n=1}^{m_2-1}
 \frac{1}{m_1^s (m_1+m_2) n^t }+\sum_{1\le m_1<m_1+m_2\le N} \frac{1}{m_1^s m_2^t(m_1+m_2) } \\
 =\sum_{m_1=1}^N \sum_{n=1}^{N-m_1-1} \sum_{n_3=1}^{N-m_1-n}
 \frac{1}{m_1^s n^t (m_1+n_3+n) } +\gz^{(N)}_\MT(s,t;1)
 \end{multline}
by setting $n_3=m_2-n$. Taking $x_1=m_1$ and $x_2=n$ in Lemma~\ref{lem:combLem}
we see that \eqref{equ:reg1} becomes
\begin{align*}
 \sum_{m_1=1}^N \sum_{n=1}^{N-m_1-1} \sum_{n_3=1}^{N-m_1-n}
& \Bigg\{\sum_{a=0}^{s-1}
{t+a-1\choose a}
 \frac{1}{m_1^{s-a} (m_1+n)^{t+a} (m_1+n_3+n) } \\
&+\sum_{b=0}^{t-1}
{s+b-1\choose b}
 \frac{1}{n^{t-b} (m_1+n)^{s+b} (m_1+n_3+n) } \Bigg\}+\gz^{(N)}_\MT(s,t;1).
\end{align*}
Taking regularization and combining with \eqref{equ:techLemmaAllgzs}
yields
\begin{align*}
 \eqref{equ:techLemma}=&
  \bgz(1,t,s)+\bgz(1,s+t)+\bgz(1,s,t)+\gz(t+1,s)-\gz_\MT(s,t;1) \\
&  -\Bigg\{\sum_{a=0}^{s-1}
{t+a-1\choose a} \bgz(1,t+a,s-a) +\sum_{b=0}^{t-1}
{s+b-1\choose b} \bgz(1,s+b,t-b)  \Bigg\}\\
=&  \bgz(1,s+t) +\gz(t+1,s)-\gz_\MT(s,t;1) \\
 & -
\Bigg\{\sum_{a=1}^{s-1}
{t+a-1\choose a} \bgz(1,t+a,s-a) +\sum_{b=1}^{t-1}
{s+b-1\choose b} \bgz(1,s+b,t-b)  \Bigg\}.
\end{align*}
Applying the stuffle relations \eqref{equ:stuffle} we get
\begin{align*}
 \eqref{equ:techLemma}=
 &\bgz(1)\gz(s+t)-\gz(s+t+1)-\gz(s+t,1)+\gz(t+1,s)-\gz_\MT(s,t;1)\\
 -& \bgz(1)\Bigg\{\sum_{a=1}^{s-1}
{t+a-1\choose a} \bgz(t+a,s-a) +\sum_{b=1}^{t-1}
{s+b-1\choose b}  \bgz(s+b,t-b)  \Bigg\}\\
+&\sum_{a=1}^{s-1}
    {t+a-1\choose a} \Big\{\gz(t+a+1,s-a)+\gz(t+a,1,s-a) \\
 &\hskip3cm   +\gz(t+a,s-a+1) +\gz(t+a,s-a,1) \Big\}\\
+&\sum_{b=1}^{t-1}
{s+b-1\choose b} \Big\{\gz(s+b+1,t-b)+ \gz(s+b,1,t-b) \\
 &\hskip3cm
+ \gz(s+b,t-b+1)+ \gz(s+b,t-b,1) \Big\}.
\end{align*}
We can cancel all the terms involving $\bgz(1)$
in the two cases $t>1$ and $t=1$
by applying \eqref{equ:EulerDec} and
\eqref{equ:EulerDec=1} of  Lemma \ref{lem:EulerDec}, respectively.
This finishes the proof of the lemma.
\end{proof}

\section{Some examples}
Using Maple we have verified all the formulas on \cite[p.~1502]{MTs}
by our general approach. We can compute any special values of
$\gz_3(s_1,\dots,s_7)$ and in particular $\gz_\slfour(s_1,\dots,s_6)$
at nonnegative integers whenever they satisfy the convergence conditions
given by \eqref{equ:convConditionOfDepth3}.
For example, by Maple computation there are 32 weight four 
convergent $\gz_\slfour(s_1,\dots,s_6)$
values at nonnegative integers, with 15 distinct ones as follows:
$$\begin{array}{ll}
 &\gz_\slfour(0,0,0,0,0,4)=\frac12\gz(2)-\frac32\gz(3)+\frac{2}{5}\gz(2)^2, \qquad 
 \gz_\slfour(0,0,0,2,2,0)=3\gz(3)-\gz(2)^2,  \\
 &\gz_\slfour(0,0,0,0,1,3)=\gz_\slfour(0,0,0,1,0,3)=\gz(2)-\gz(3)-\frac{1}{10}\gz(2)^2, \\ 
& \gz_\slfour(0,0,0,0,2,2)=\gz_\slfour(0,0,0,2,0,2)=\gz(3)-\frac{3}{10}\gz(2)^2, \\
 &\gz_\slfour(1,0,0,0,1,2)=\gz_\slfour(0,0,1,1,0,2)=\gz(3)-\frac{1}{5}\gz(2)^2,\\
 &\gz_\slfour(1,0,0,0,2,1)=\gz_\slfour(0,0,1,2,0,1)=2\gz(3)-\frac{1}{2}\gz(2)^2,\\
& \gz_\slfour(0,0,1,0,0,3)=\gz_\slfour(0,1,0,0,0,3)=\gz_\slfour(1,0,0,0,0,3)
=2\gz(3)-\gz(2)-\frac{1}{10}\gz(2)^2, \\ 
& \gz_\slfour(0,1,0,1,1,1)=\frac{7}{10}\gz(2)^2,\qquad \gz_\slfour(1,1,1,0,0,1)=\frac{12}{5}\gz(2)^2,\\
& \gz_\slfour(0,0,0,1,1,2)=\frac{1}{10}\gz(2)^2,\qquad \gz_\slfour(1,0,1,1,1,0)=\frac{17}{10}\gz(2)^2,\\
 &\gz_\slfour(1,0,1,0,0,2)=\gz_\slfour(0,1,1,0,0,2)=\gz_\slfour(1,1,0,0,0,2)=\frac{4}{5}\gz(2)^2,\\
 &\gz_\slfour(1,0,0,1,1,1)=\gz_\slfour(1,0,0,1,2,0)=\gz_\slfour(0,0,1,2,1,0)=\gz_\slfour(0,0,1,1,1,1)=\frac{1}{2}\gz(2)^2,\\
 &\gz_\slfour(0,1,0,1,0,2)=\gz_\slfour(1,0,0,1,0,2)=\gz_\slfour(0,1,0,0,1,2)=\gz_\slfour(0,0,1,0,1,2)=\frac{2}{5}\gz(2)^2,\\
 &\gz_\slfour(1,1,0,0,1,1)=\gz_\slfour(1,0,1,1,0,1)=\gz_\slfour(1,0,1,0,1,1)=\gz_\slfour(0,1,1,1,0,1)=\frac{6}{5}\gz(2)^2.
\end{array} $$
Notice that in the first seven values the weights are mixed which correspond 
to the nine cases of irregular values \eqref{equ:irr1} to \eqref{equ:irr5}. 
In higher weight cases if we only consider regular (pure weight) values then
we have the following examples.
$$\begin{array}{ll}
&\gz_\slfour(1,1,0,1,1,1)=\gz_\slfour(0,1,1,1,1,1)=\frac52\gz(5)-\gz(2)\gz(3)=.6150150376
 \dots,\\
& \gz_\slfour(1,0,1,1,1,1)=-\frac32\gz(5)+\gz(2)\gz(3)=.4219127176 \dots,\\
& \gz_\slfour(1,1,1,1,1,1)= -\frac{62}{105}\gz(2)^3+2\gz(3)^2=.2617453537 \dots,\\
& \gz_3(1,1,1,1,1,1,1)=\frac32\gz_\slfour(1,1,1,1,1,2)= \frac{21}8\gz(7)-\frac32\gz(2)\gz(5)=.08840016918 \dots,\\
& \gz_\slfour(2,2,2,2,2,2)=\frac{368}{875875}\gz(2)^6=.0083233212. 
\end{array} $$ 
The last value of the above agrees with \cite[(4.29)]{KMT0} and
the first value in \cite[Table 4]{GS} by the relation
$$\gz_\slfour( 2m,2m,2m,2m,2m,2m )=12^{-2m} \gz_\slfour(2m),$$
where $\gz_\slfour$ on the LHS is the multiple variable function 
being studied in the current paper
and $\gz_\slfour$ on the RHS is the original single variable 
Witten zeta function
(see \cite[Prop.~2.1]{MTs}). As two intriguing computation we have
$$\begin{array}{ll}
&\gz_\slfour(1,2,3,3,2,1)=\gz(2)\gz(8,2)+\frac{811324}{238875}\gz(2)^6
    -\frac{5}{2}\gz(2)\gz(5)^2-\frac{37}{2}\gz(3)\gz(9)-35\gz(5)\gz(7)\\
    \ & \hskip6.5cm-2\gz(7)\gz(2)\gz(3)+\frac{37}{4}\gz(10,2)=.0129650292\dots, \\
& \gz_\slfour(3,2,1,1,2,3)=10\gz(2)\gz(8,2)-\frac{120112}{53625}\gz(2)^6
    -6\gz(2)\gz(5)^2+44\gz(3)\gz(9) +40\gz(5)\gz(7)\\
   \ & \hskip6.5cm -20\gz(7)\gz(2)\gz(3)-22\gz(10,2)=.0056078053\dots.
\end{array} $$
We have also verified all the above computation
numerically with Maple and EZface \cite{EZface}. 

We conclude our paper with a conjecture.
Let $\MZV(w,\le l)$ be the $\Q$-vector space generated by MZVs of weight  $w$
and depth $\le l$. Then by Theorem~\ref{thm:depth3} the space generated by
special values of $\gz_\slfour$ of weight $w>3$ over $\Q$ is included 
in the space $\MZV(w,\le 3)+\MZV(w-1,\le 2)+\MZV(w-2,1)$ (conjecturally the sum is direct). 
\begin{conj} The space generated by
special values of $\gz_\slfour$ of weight $w>3$ at nonnegative integers
over $\Q$ is 
$$\MZV(w,\le 3)\oplus \MZV(w-1,\le 2)\oplus  \MZV(w-2,1).$$
\end{conj}
Using the MZV table \cite{pet} we have verified the conjecture for all weights up to 12.

\end{document}